\documentstyle[amssymb,12pt]{amsart}
\numberwithin{equation}{subsection}
 
\newcommand{\rk}{\operatorname{rk}}
\newcommand{\kap}{\kappa}
\newcommand{\th}{\theta}

\newcommand{\ga}{\gamma}
\newcommand{\varth}{\vartheta}

\renewcommand{\H}{{\goth H}}
\newcommand{\SL}{\operatorname{SL}}
\newcommand{\GL}{\operatorname{GL}}

\renewcommand{\Im}{\operatorname{Im}}
\newcommand{\ra}{\rightarrow}
\newcommand{\EE}{{\cal E}}
\newcommand{\VV}{{\cal V}}

\newcommand{\PP}{{\cal P}}
\newcommand{\A}{{\cal A}}
\newcommand{\LL}{{\cal L}}
\renewcommand{\O}{{\cal O}}

\newcommand{\om}{\omega}
\newcommand{\la}{\lambda}
\newcommand{\C}{{\Bbb C}}

\newcommand{\Z}{{\Bbb Z}}

\newcommand{\Ga}{\Gamma}

\newcommand{\wt}{\widetilde}
\newcommand{\ot}{\otimes}

\newcommand{\sub}{\subset}
 
\title{M.~P.~Appell's function and vector bundles of rank 2 on elliptic curves}
\begin{document}
\maketitle
\centerline{A. Polishchuk\footnote{E-mail: apolish@@math.harvard.edu}}

\bigskip

\centerline{Department of Mathematics}
\centerline{Harvard University, Cambridge, MA 02138}       

\vspace{6mm}

The purpose of this work is to establish a connection between
the function
$$
\kap(y,x)=\kap(y,x,\tau)=\sum_{n\in\Z}\frac{\exp(\pi i\tau n^2+
2\pi i n x)}{\exp(2\pi i n\tau)-\exp(2\pi i y)},
$$
where $y,x,\tau\in\C$, $\Im(\tau)>0$, $y\not\in \Z+\Z\tau$,
and vector
bundles of rank 2 on elliptic curves. This function was
introduced by M.~P.~Appell
in \cite{Ap} in order to decompose into simple elements
the so called elliptic functions of the third
kind (which correspond to meromorphic sections of line bundles
on elliptic curves).\footnote{The same function
was introduced by M.~Hermite essentially for the same problem,
however, he didn't publish his results until the appearance
of the first part of Appell's paper. In the second part
of this paper Hermite's method is partially reproduced.
Appell actually used the function which differs from $\kap$ by a
theta-function. Our notation is closer to that of Halphen's
book \cite{H}.}
 
As is well-known Riemann's theta function arises naturally when
considering
global  sections  of  line  bundles on elliptic curves.
We claim that the function  $\kap$  is  connected  in  a
similar  way  with  rank-2  bundles.  Namely, the difference
equation for $\kap$ allows  to  interpret  it  as  a  global
section of a rank-2 bundle of degree 1 on an elliptic curve.
Pursuing this analogy we derive some  interesting
identities satisfied  by  $\kap$  similar  to  the  addition
formulas for $\th$. These identities (that were known to Appell) 
appear as $A_{\infty}$-constraints
in the Fukaya category of an elliptic curve (see \cite{P}).
Another identity for function $\kap$ turned out to be
useful in proving some formulas from Ramanujan's Notebooks
(cf. \cite{E}). 
However, at the present moment there appears to be no general
theory similar to that of theta-identities. 
The following pair of identities is derived easily from those 
known to Appell:
$$\th(0)\kap(\tau/2,1/2)+\th(1/2)\kap((\tau+1)/2,0)=
\frac{1}{2}
\th(\tau/2)^3.$$
$$\th(\tau/2)^3\kap(1/2,\tau/2)=\th(1/2)^3\kap(\tau/2,1/2)+
\th(0)^3\kap((\tau+1)/2,0).$$
These identities relating three special values $\kap(y,x)$
with $y,x$ points of order 2 on elliptic curve $\C/\Z+\tau\Z$ deserve to be 
considered as analogues of
the famous Jacobi identity relating the 4-th powers of values of
$\theta$ at points of order 2. In fact, all the values 
$\kap(x,y)$, where $y$ and $x$ are points of order 2 on elliptic
curve, except these three can be expressed in terms of theta-function.

One may ask whether there is an analogue
of modular property for the function $\kap$ similar to the
functional equation for the theta function. Note that the
latter equation reflects the fact that one can construct a
line bundle  on  the  universal  elliptic  curve  (over  the
analytic moduli stack of elliptic curves with a non-trivial
point of order 2) such that $\th$ is its section. Similar,
we can construct a vector bundle of rank 2 on the universal
elliptic curve. However, it seems that there is no simple
formula for the corresponding 1-cocycle of the modular group.
So  we  propose  to  formulate  the  modular property of the
function $\kap$ not in the form of the equation but  in  the
form of divisibility property. More precisely, let us denote
$$\kap_0(x,\tau)=\exp(\frac{3\pi i\tau}{4})
\kap((\tau+1)/2,x,\tau).$$
We prove that for every element
$\ga=
\left(\matrix a & b \\ c & d \endmatrix\right)$
of the group $\Ga_{1,2}$
(this is the standard subgroup of $\SL_2(\Z)$ which preserves
the point $(\tau+1)/2$) the difference
$$\kap_0(\frac{x}{c\tau+d},\frac{a\tau+b}{c\tau+d})-
\zeta\cdot (c\tau+d)\cdot\exp(\pi i(\frac{1}{c\tau+d}-1)x)
\kap_0(x,\tau)$$
is divisible by $\th(x,\tau)$ in  the  ring  of  holomorphic
functions  on $\H\times\C$, where $\zeta$ is a root of unity
of order 4 which we explicitly determine.

Although we concentrate on the case of bundles of degree 1,                
it seems that almost all facts concerning vector bundles of rank 2
on elliptic curves can be worked out explicitly in terms of the
function $\kap$ (and elliptic functions). 
It seems that the case of higher rank bundles is described
similarly by various functions derived from $\kap$
by taking derivatives and difference
derivatives $(\kap_{a_1}-\kap_{a_2})/(a_1-a_2)$, etc.

In section \ref{isom} we show that some natural isomorphisms between
rank-2 vector bundles can be written explicitly in terms of $\kap$.
As a byproduct we show how to find explicitly (in terms of $\kap$)
holomorphic functions $\phi_1$ and
$\phi_2$ on $\C^*$ such that
$$\phi_1(z)\th(z,q^2)-\phi_2(z)\th(qz,q^2)=1$$
(one  knows  apriori  that  such   functions   exist   since
$\th(z,q^2)$ and $\th(qz,q^2)$ have no common zeroes).

We use both  additive  and  multiplicative  notations  (i.e.
sometimes  we  consider  variables  in $\C$ and sometimes in
$\C^*$). The relation between (some) multiplicative and additive
variables is the following: $z=\exp(2\pi i x)$,
$q^{1/2}=\exp(\pi i \tau)$, $a=\exp(2\pi iy)$. By abuse of
notation, for any function in multiplicative variables we will denote
the corresponding function in additive variables by the same
letter (the only section where the additive notation is
used is section \ref{modular} on modular property of $\kap$).
When working with a fixed elliptic curve $E_q=\C^*/q^{\Z}$ we often omit
the variable $q$ in the corresponding elliptic functions and
the function $\kap$.

\noindent {\it Acknowledgment}. I'd like to thank B.~Berndt and
G.~Andrews for help with finding the references.
This work is partially supported by NSF grant.
 
\section{Bundles of rank 2 and degree 1}\label{rank2}
        
In this section we always use multiplicative variables.
 
\subsection{}
Let $q\in\C^*$ be such that $|q|<1$,
$E=E_q=\C^*/q^{\Z}$ be the corresponding elliptic curve.
For every invertible $r\times r$ matrix $A(z)$ of holomorphic
functions  on  $\C^*$  we  define  the  corresponding
holomorphic vector bundle $V_r(A)$ on $E_q$ as follows:
$$V_r(A)=\C^*\times\C^r/(z,v)\sim(qz,A(z)v).$$
Holomorphic sections of $V_r(A)$ correspond to
$r$-tuple of functions $v(z)$ satisfying the equation
$$v(qz)=A(z)v(z).$$
The bundles $V_r(A)$ and $V_r(A')$ are isomorphic if and only if
there exists an invertible $r\times r$ matrix $B(z)$ such
that $A'(z)=B(qz)A(z)B(z)^{-1}$. The construction $A\mapsto
V(A)$ is compatible with tensor products, in particular,
if $\phi$ is an invertible holomorphic function on $\C^*$
then
$$V_r(\phi A)\simeq V_1(\phi)\otimes V_r(A).$$
        
\subsection{} We denote by $L$   the   line   bundle
$V_1(q^{-1/2}z^{-1})$. The classical theta function
$$\th(z)=\th(z,q)=\sum_{n\in\Z}q^{n^2/2}z^n$$
is a section of $L$ (as function in $z$). We also denote by $P_a$
the line bundle $V_1(a)$, where $a\in\C^*$ is considered as a constant
function on $\C^*$.
Then we have a canonical isomorphism
\begin{equation}\label{sq1}
P_a\otimes t^*_aL\simeq L.
\end{equation}
In particular, this implies that for every $a$ one has
\begin{equation}\label{add}
t_a^*L\otimes t_{a^{-1}}^*L\simeq L^2.
\end{equation}
The classical addition formula for theta-function which we
remind below
gives a more concrete form of this isomorphism.
The basis of global sections of $L^2$ consists of
functions $\varth_0$ and $\varth_1$ defined as follows:
$$\varth_0(z)=\varth_0(z,q)=\th(z^2,q^2)=
\sum_{n\in\Z}q^{n^2}z^{2n},$$
$$\varth_1(z)=\varth_1(z,q)=\th
\left[\matrix 1/2 \\ 0 \endmatrix\right](z^2,q^2)=
\sum_{n\in\Z}q^{(n+1/2)^2}z^{2n+1}.$$
Now the addition formula corresponding to (\ref{add}) is
\begin{equation}\label{addf}
\th(za)\th(za^{-1})=
\varth_0(a)\varth_0(z)+\varth_1(a)\varth_1(z).
\end{equation}
       
\subsection{} Consider the rank-2 bundle
$$F_a=V_2
\left(  \matrix  a  &  1  \\  0  & q^{-1/2}z^{-1} \endmatrix
\right).$$
This bundle is a unique non-trivial extension of $L$
by  $P_a$. Assume that $a\not\in q^{\Z}$.
Then $H^0(E,P_a)=H^1(E,P_a)=0$, hence
the natural projection $H^0(E,F_a)\ra H^0(E,L)$ is an isomorphism.
This means that there exists a unique global section of $F_a$
projecting   to  $\th(z,q)$.  This section has form
$$\left[\matrix \kap_a(z) \\ \th(z) \endmatrix\right]$$
where $\kap_a(z)=\kap(a,z)=\kap(a,z,q)$
is the unique holomorphic in $z$ function satisfying
the equation
\begin{equation}\label{def}
\kap(a,qz,q)=a\kap(a,z,q)+\th(z,q).
\end{equation}
This function can be written as the following series
$$\kap(a,z,q)=\sum_{n\in\Z}\frac{q^{n^2/2}}{q^n-a}z^n$$
converging for all $z$, since $|q|<1$.
Using the defining equation (\ref{def}) one can easily
establish the following identities:
\begin{equation}\label{inv}
\kap(a,z)=-a^{-1}\kap(a^{-1},qz^{-1}),
\end{equation}
\begin{equation}\label{def2}
\kap(qa,z)=q^{-1/2}z\kap(a,qz)=q^{-1/2}az\kap(a,z)
+q^{-1/2}z\th(z)
\end{equation}
The latter one allows to consider the pair $(\kap,\th)$
as a meromorphic section
of  certain  rank-2  bundle  on  $E\times  E$.  It  is  also
convenient to consider the function
$$\bar{\kap}(a,z)=\th(-q^{-1/2}a)\kap(a,z)$$
which is holomorphic in $a$ and $z$. Comparing the quasi-periodicity
properties of this function in $a$ and $z$ one immediately
arrives to the following symmetry identity:
\begin{equation}\label{sym}
a\bar{\kap}(a,z)=-q^{1/2}z\bar{\kap}(-q^{1/2}z,-q^{-1/2}a).
\end{equation}
Note that this identity appears as the main theorem of \cite{E}
where it is applied to derive some of Ramanujan's formulas.
However, it can be found already in Appell's paper \cite{Ap}
(t.III, p.20) where it is shown to be equivalent to certain identity
of Hermite.
 
Finally, the function $\bar{\kap}(az,z^{-1})$ for fixed $a$
is a section of $P_{-q^{-1/2}a^{-1}}\otimes L^2$. If we choose
a square root $a^{1/2}$ we can write this as
$t^*_{iq^{1/4}a^{1/2}}L^2$. Hence, we can express
$\bar{\kap}(az,z^{-1})$ as a linear combination
of the corresponding shifts of $\varth_0$ and $\varth_1$.
Using the vanishing $\varth_0(iq^{1/2})=\varth_1(i)=0$
we  find  explicitly  the coefficients and get the following
identity:
\begin{equation}\label{sqrt}
\bar{\kap}(az,z^{-1})=\frac{\varth_0(iq^{1/4}a^{1/2}z)}
{\varth_0(i)}\bar{\kap}(q^{-1/4}a^{1/2},q^{1/4}a^{1/2})
+\frac{\varth_1(iq^{1/4}a^{1/2}z)}
{\varth_1(iq^{1/2})}\bar{\kap}(q^{1/4}a^{1/2},q^{-1/4}a^{1/2})
\end{equation}
       
\subsection{} The analogue of the property (\ref{sq1})
for rank-2 bundles $F_a$ is the following canonical
isomorphism:
\begin{equation}\label{sq2}
t_b^*F_a\simeq P_{b^{-1}}\otimes F_{ab}.
\end{equation}
Indeed, the left hand side is
$$V
\left( \matrix a & 1\\ 0 &  b^{-1}q^{-1/2}z^{-1}  \endmatrix
\right)$$
while the right hand side is
$$V
\left( \matrix a & b^{-1}\\ 0 &  b^{-1}q^{-1/2}z^{-1}  \endmatrix
\right)$$
Now the above isomorphism follows from the equality
$$
\left( \matrix 1 & 0\\ 0 &  b \endmatrix\right)
\left( \matrix a & 1\\ 0 &  b^{-1}q^{-1/2}z^{-1}  \endmatrix
\right) \left( \matrix 1 & 0\\ 0 &  b^{-1} \endmatrix\right)
=\left( \matrix a & b^{-1}\\ 0 &  b^{-1}q^{-1/2}z^{-1}  \endmatrix
\right).$$
Using (\ref{sq1}) and (\ref{sq2}) we obtain an isomorphism
\begin{equation}\label{sq3}
L\ot t_b^*F_a\wt{\ra} t_b^*L\ot F_{ab}.
\end{equation}
Moreover,  this  isomorphism is compatible with morphisms of
both bundles to $t_b^*L\ot L$. Now the section
$$\left[\matrix \th\cdot t_{b}^*\kap_a \\
\th\cdot t_b^*\th \endmatrix\right]$$
of $L\ot t_b^*F_a$ is
mapped by (\ref{sq3}) to the section
$$\left[\matrix \th\cdot t_{b}^*\kap_a \\
b^{-1}\th\cdot t_b^*\th \endmatrix\right]$$
of $t_b^*L\ot F_{ab}$.
On the other hand, the latter bundle has the global section
$$\left[\matrix \kap_{ab}\cdot t_b^*\th \\
\th\cdot t_b^*\th \endmatrix\right].$$
It follows that the difference
$$\kap_{ab}\cdot t_b^*\th - b^{-1}\th\cdot t_b^*\kap_a$$
is a constant multiple of $t_{a^{-1}}^*\th$.
Considering the value of this difference at $z=-q^{1/2}$
and using vanishing $\th(-q^{1/2})=0$ we get the following
identity
\begin{equation}\label{h-add}
\th(bz)\kap(ab,z)-b^{-1}\th(z)\kap(a,bz)=\frac{\th(-q^{1/2}b)
\th(a^{-1}z)}{\th(-q^{-1/2}a)}\kap(ab,-q^{1/2}).
\end{equation}
This identity appears in  a  slightly  different  form in
Halphen's  book  \cite{H}  (p.481, formula (45) and the next
one).
Note that $\bar{\kap}(a,-q^{1/2})=a\bar{\kap}(a,-q^{-1/2})$
and (\ref{sym}) gives
$$a\bar{\kap}(a,-q^{-1/2})=\bar{\kap}(1,-q^{-1/2}a).$$
Now using the explicit series for  $\kap(a,z)$  we
immediately see that
$$\bar{\kap}(1,z)=\lim_{a\ra 1}\frac{\th(-q^{-1/2}a)}{1-a}=
q^{-1/2}\frac{d\th}{da}(-q^{-1/2})=
\frac{1}{2}\th(1)\th(-1)\th(q^{1/2})$$
due to Jacobi's derivative formula.
It follows that
\begin{equation}\label{sp1}
\kap(a,-q^{1/2})=\frac{\th(1)\th(-1)\th(q^{1/2})}
{2\th(-q^{-1/2}a)}.
\end{equation}
Substituting this into (\ref{h-add}) we get
\begin{equation}\label{h-add2}
\th(bz)\kap(ab,z)-b^{-1}\th(z)\kap(a,bz)=\frac{
\th(1)\th(-1)\th(q^{1/2})\th(-q^{1/2}b)\th(a^{-1}z)}
{2\th(-q^{-1/2}a)\th(-q^{-1/2}ab)}.
\end{equation}
This identity can be considered as a rank-2 analogue of the addition
formula (\ref{addf}).
It allows to express $\kap$ in terms of $\kap(q^{1/2},z)$ and
theta-functions:
\begin{equation}\label{h-add3}
\kap(a,z)=q^{1/2}a^{-1}\frac{\th(z)}{\th(q^{-1/2}az)}
\kap(q^{1/2},q^{-1/2}az)+\frac{\th(1)\th(q^{1/2})
\th(-a)\th(q^{-1/2}z)}{2\th(-q^{-1/2}a)\th(q^{-1/2}az)}.
\end{equation}
  
\subsection{} We have a canonical morphism
$$\mu_{a,b}:H^0(E,t_{b^{-1}}^*L)\otimes H^0(E,t_b^*F_a)
\ra H^0(E, t_{b^{-1}}^*L\ot t_b^*F_a)\simeq H^0(E, L\otimes F_{ab}).$$
Assume that $a$ and $ab$ are not integer powers of $q$.
Then this map is given by
$$\mu_{a,b}(t_{b^{-1}}^*\th\otimes
\left[\matrix t_{b}^*\kap_a \\ t_b^*\th \endmatrix\right])=
\left[\matrix t_{b^{-1}}^*\th\cdot t_{b}^*\kap_a \\
b t_{b^{-1}}^*\th \cdot t_b^*\th \endmatrix\right].$$
Using the maps $\mu_{a,b}$ one can easily construct
a basis in $H^0(E,L\otimes F_a)$.
Indeed, assume that $a\in\C^*$ is such that $a\not\in q^{\Z}$
and $-a\not\in q^{\Z}$. Then we claim that
$H^0(E, L\otimes F_a)$ is generated
by the images of $\mu_{a,1}$, $\mu_{-a,-1}$ and the subspace
$H^0(E,L\otimes P_a)$.
This follows immediately from the exact sequence
\begin{equation}\label{exact}
0\ra H^0(E, L\otimes P_a)\ra H^0(E, L\otimes F_a)\ra
H^0(E, L^2)\ra 0.
\end{equation}
and the fact that $\th^2$ and $t_{-1}^*\th^2$ generate $H^0(E,L^2)$.
It follows that the following vectors
form a basis of $H^0(E, L\otimes F_a)$:
$$v_0(a)=\left[\matrix \th(a^{-1}z) \\
0 \endmatrix\right],
v_1(a)=\left[\matrix \th(z)\cdot \kap(a,z) \\
\th(z)^2 \endmatrix\right],
v_{-1}(a)=\left[\matrix \th(-z)\cdot \kap(-a,-z) \\
-\th(-z)^2 \endmatrix\right].$$
  
One can write the map $\mu_{a,b}$ expicitly in terms of this basis
using (\ref{h-add3}). Namely, one has
$$b^{-1}\cdot\mu_{a,b}(t_{b^{-1}}^*\th\otimes
\left[\matrix t_{b}^*\kap_a \\ t_b^*\th \endmatrix\right])=
\la_b\cdot v_1(ab)-\la_{-b}\cdot v_{-1}(ab)+(\nu_{a,b}-\nu_{a,-b})\cdot
v_0(ab)$$
where
$$\la_b=\frac{\th(q^{1/2}b^{-1})\th(q^{1/2}b)}{\th(q^{1/2})^2},$$
$$\nu_{a,b}=\frac{\th(1)\th(q^{1/2}b)\th(-q^{1/2}b)\th(-q^{1/2}b^{-1})
\th(b)\th(ab)}{2\th(q^{1/2})\th(-q^{1/2}a^{-1})\th(q^{-1/2}ab)\th(-ab^2)}.$$
  
\section{Various formulas}
 
\subsection{} Let us consider special values of
$\bar{\kap}(a,z)$ for $a,z\in\{\pm 1,\pm q^{1/2}\}$.
It is easy to see that all of them except for
$\bar{\kap}(q^{1/2},-1)$, $\bar{\kap}(-q^{1/2},1)$ and
$\bar{\kap}(-1,q^{1/2})$
can be expressed in terms of values of the theta function.
We have already seen how to compute $\kap(a,-q^{1/2})$. For example,
\begin{equation}\label{sp2}
\kap(-q^{1/2},-q^{1/2})=\frac{1}{2}\th(-1)\th(q^{1/2}),
\end{equation}
\begin{equation}\label{sp3}
\kap(-1,-q^{1/2})=\frac{1}{2}\th(1)\th(-1).
\end{equation}
The latter identity appears as entry 18 in \cite{B}, p.152.
Also using (\ref{inv}) and (\ref{sym}) one can easily deduce
that
\begin{equation}\label{sp4}
\kap(-1,\pm 1)=\frac{1}{2}\th(\pm 1),
\end{equation}
\begin{equation}\label{sp5}
\kap(\pm q^{1/2}, q^{1/2})=\frac{1}{2}\th(q^{1/2}),
\end{equation}
while $\kap(q^{1/2},1)=\kap(-q^{1/2},-1)=0$ as one can see
immediately from the formulas
$$\kap(q^{1/2},z)=\sum_{n\ge 0}\frac{q^{n^2/2+n}}
{1-q^{n+1/2}}(z^{-n}-z^{n+1}),$$
$$\kap(-q^{1/2},z)=\sum_{n\ge 0}\frac{q^{n^2/2+n}}
{1+q^{n+1/2}}(z^{-n}+z^{n+1}).$$
Note that three exceptional values correspond to
pairs $(a,z)$ which are stable under involution
$(a,z)\mapsto (-q^{1/2}z, -q^{-1/2}a)$ (modulo $q^{\Z}$).
 
\subsection{} Using the above information about the  special
values of $\kap$ one can derive from (\ref{h-add3}) the
following identities
\begin{equation}\label{id4}
2\kap(q^{1/2}b^{-1},q^{1/2}b)=\th(q^{-1/2}b)+
\frac{\th(1)\th(b)}{\th(-b)}\th(-q^{1/2}b).
\end{equation}
\begin{equation}\label{id5}
\kap(q^{1/2}b^{-1},b)=\frac{\th(q^{1/2})\th(q^{-1/2}b)
\th(-q^{-1/2}b)}{2\th(-b)}=
\frac{(q)_{\infty}^2(-q)^2_{\infty}(-b)_{\infty}(-qb^{-1})_{\infty}
(b)_{\infty}(qb^{-1})_{\infty}}{(q^{1/2}b)_{\infty}
(q^{1/2}b^{-1})_{\infty}}.
\end{equation}
\begin{equation}\label{id5.5}
\th(-z)\kap(a,z)+\th(z)\kap(-a,-z)=
\frac{\th(q^{1/2})^2\th(1)\th(-1)\th(-a^{-1}z)}{2\th(q^{1/2}a^{-1})
\th(-q^{1/2}a^{-1})}.
\end{equation}
Here is another identity which follows from (\ref{h-add3}):
\begin{equation}\label{id6}
\th(q^{1/2})^2\th(-b)\kap(q^{1/2}b^{-1},-b)=
\th(q^{1/2}b^{-1})^2\th(-1)\kap(q^{1/2},-1)+
\th(-q^{1/2}b^{-1})^2\th(1)\kap(-q^{1/2},1)
\end{equation}
 
\subsection{} The special values $\kap(q^{1/2},-1)$, $\kap(-q^{1/2},1)$
and $\kap(-1,q^{1/2})$, which we were not able to express in terms
of theta-functions, satisfy the following two identities: 
\begin{equation}\label{for1}
\th(1)\kap(q^{1/2},-1)+\th(-1)\kap(-q^{1/2},1)=\frac{1}{2}
\th(q^{1/2})^3.
\end{equation}
\begin{equation}\label{for2}
\th(q^{1/2})^3\kap(-1,q^{1/2})=\th(-1)^3\kap(q^{1/2},-1)+
\th(1)^3\kap(-q^{1/2},1).
\end{equation}
The first formula is obtained by substituting
$b=z=-1$, $a=-q^{1/2}$ into (\ref{h-add2}), while (\ref{for2}) 
is the specialization of (\ref{id6}) for $b=-q^{1/2}$.
We can rewrite (\ref{for1}) as follows
\begin{align*}
&(\sum_{k\in\Z}q^{k^2/2})\left(\sum_{n\ge0}(-1)^n
\frac{q^{n^2/2+n}}{1-q^{n+1/2}}\right)+
(\sum_{k\in\Z}(-1)^kq^{k^2/2})\left(\sum_{n\ge0}
\frac{q^{n^2+n}}{1+q^{n+1/2}}\right)=\\
&2\left(\sum_{n\ge0}q^{(n^2+n)/2}\right)^3.
\end{align*}
The LHS can be further rewritten as 
\begin{align*}
&\sum_{n\ge 0,k\in n+2\Z}(-1)^nq^{k^2/2+n^2/2+n}(\frac{1}{1-q^{n+1/2}}
+\frac{1}{1+q^{n+1/2}})+\\
&\sum_{n\ge 0,k\in n+1+2\Z}(-1)^nq^{k^2/2+n^2/2+n}(\frac{1}{1-q^{n+1/2}}
-\frac{1}{1+q^{n+1/2}})=\\
&2\sum_{n\ge 0,l\in\Z}(-1)^n\frac{q^{(n+2l)^2/2+n^2/2+n}}{1-q^{2n+1}}
+2\sum_{n\ge 0,l\in\Z}(-1)^n\frac{q^{n-2l+1)^2/2+n^2/2+2n+1/2}}{1-q^{2n+1}}.
\end{align*}
Simplifying we get the identity
$$\left(\sum_{n\ge0}q^{(n^2+n)/2}\right)^3=\sum_{n\ge 0,l\in\Z}
(-1)^n\frac{q^{n^2+n-2nl+2l^2}(1+q^{2l+1})}{1-q^{2n+1}}$$
Note that G.~Andrews in \cite{A} has obtained another formula for
$(\sum q^{(n^2+n)/2})^3$, namely
$$\left(\sum_{n\ge0}q^{(n^2+n)/2}\right)^3=\sum_{n\ge 0,2n\ge j\ge0}
\frac{q^{2n^2+2n-j(j+1)/2}(1+q^{2n+1})}{1-q^{2n+1}}$$
where all the coefficients in the RHS are positive (which in particular
proves that every number is a sum of three triangular numbers).
It would be interesting to see directly why the right hand sides in
two formulas are equal.
                                    
\section{Isomorphisms between bundles}\label{isom}

In this section we show that some isomorphisms between
rank-2 bundles can be written explicitly using the function
$\kap$.

\subsection{} The first isomorphism we are interested in
is an isomorphism of $F_a$ (see section \ref{rank2}) and
the unique non-trivial extension of $\det F_a\simeq P_a\ot L$
by $\O$. The latter extension is realized by the rank-2 bundle
$$F'_a=V_2
\left(  \matrix  1 &  1  \\  0  & q^{-1/2}az^{-1} \endmatrix
\right).$$
To give an isomorphism $F'_a\simeq F_a$ we have to
find a $2\times 2$ matrix of holomorphic functions $B(z)$ with
the property
\begin{equation}\label{conj1}
\left(  \matrix  a  &  1  \\  0  & q^{-1/2}z^{-1} \endmatrix
\right)=B(qz)
\left(  \matrix  1 &  1  \\  0  & q^{-1/2}az^{-1} \endmatrix
\right)B(z)^{-1}.
\end{equation}
Since $B$ should send global sections of $F'_a$ to global
sections of $F_a$ we can choose $B$ in the form
$$B=\left(  \matrix  \kap_a(z) &  f(z)  \\  \th(z) & g(z) \endmatrix
\right).$$
Now it is easy to see that $(\ref{conj1})$ implies
$$g(qz)=a^{-1}g(z)-a^{-1}\th(z),$$
hence, $g(z)=-a^{-1}\kap_{a^{-1}}(z)$.
To find $f(z)$ we note that the determinant of $B$ should be
$q$-periodic, hence, a constant function. Therefore, we have
$$f(z)\th(z)=c-a^{-1}\kap_a(z)\kap_{a^{-1}}(z)$$
for some constant $c$. Substituting $z=-q^{1/2}$ we
find using (\ref{sp1}) that
$$c=c_a=a^{-1}\kap_a(-q^{1/2})\kap_{a^{-1}}(-q^{1/2})=
\frac{\th(1)^2\th(-1)^2\th(q^{1/2})^2}{4a\th(-q^{-1/2}a)
\th(-q^{1/2}a)},$$
and the matrix $B$ is equal to
$$\left(  \matrix \kap_a(z)  &  \frac{c_a-a^{-1}\kap_a(z)
\kap_{a^{-1}}(z)}{\th(z)}\\  \th(z)  & -a^{-1}
\kap_{a^{-1}}(z) \endmatrix \right).$$
                            
\subsection{} Now we want to compare bundles $F_a\simeq F'_a$
with push-forwards of line bundles of degree 1 on
$E_{q^2}=\C^*/q^{2\Z}$ under the natural isogeny
$\pi:E_{q^2}\ra E_q$. It is sufficient to do this for $F'_1$
since $F'_a$ is obtained from $F'_1$ by translation.
We claim that there is an isomorphism
$$F'_1\simeq\pi_*L'$$
where $L'$ is a line bundle on $E_{q^2}$ with multiplicator
$-q^{-1/2}z^{-1}$, i.e.
$$L'=V^{q^2}_1(-q^{-1/2}z^{-1}),$$
where $V^{q^2}_1$ is the construction of section \ref{rank2}
with $q$ replaced by $q^2$. Indeed, the unique section of $L'$
(which is given by $\th(-q^{-1/2}z,q^2)$) induces the global
non-vanishing section of $\pi_*L'$, so $\pi_*L'$ fits into
exact sequence
$$0\ra \O\ra\pi_*L'\ra L\ra 0.$$
On  the  other  hand,  it is easy to check that $\pi_*L'$ is a
simple bundle, hence the above sequence doesn't split and $\pi_*L'$
is isomorphic to $F'_1$. To write this isomorphism explicitly
note that
$$\pi_*L'=V_2
\left(  \matrix  0 & -q^{-1/2}z^{-1} \\  1  & 0 \endmatrix
\right).$$
Thus, we are looking for a matrix $C(z)$ such that
\begin{equation}\label{conj2}
\left(  \matrix  1  &  1  \\  0  & q^{-1/2}z^{-1} \endmatrix
\right)=C(qz)
\left(  \matrix  0 &  -q^{-1/2}z^{-1} \\  1  & 0 \endmatrix
\right)C(z)^{-1}.
\end{equation}
Let $C(z)=\left(\matrix c_{11}(z) & c_{12}(z) \\  c_{21}(z) &
c_{22}(z)  \endmatrix\right).$  Then  from  (\ref{conj2}) we
deduce that the vector
$$\left[\matrix c_{22}(z) \\ -c_{21}(z)\endmatrix\right]$$
is a global section of $\pi_*L'$. Hence, we can choose $C$
with $c_{21}(z)=\th(-q^{-1/2}z,q^2)$ and
$c_{22}(z)=-c_{21}(qz)$. As before $(\ref{conj2})$
implies that $\det(C)$ is a (non-zero) constant, i.e.
\begin{equation}\label{deteq}
c_{11}(z)c_{21}(qz)+c_{12}(z)c_{21}(z)=c
\end{equation}
for some constant $c$. Another consequence of (\ref{conj2})
is that $c_{11}(z)=c_{12}(qz)-c_{21}(z)$. Substituting this
in (\ref{deteq}) we get
$$c_{12}(qz)c_{21}(qz)=-c_{12}(z)c_{21}(z)+c+
c_{21}(z)c_{21}(qz).$$
In other words, if we denote $\phi(z)=c_{12}(z)c_{21}(z)-c/2$
then we have
$$\phi(qz)=-\phi(z)+c_{21}(z)c_{21}(qz).$$
Note that
$$c_{21}(z)c_{21}(qz)=\la\cdot\th(-z,q)$$
where
$$\la=\frac{\th(1,q^2)\th(q,q^2)}{\th(q^{1/2},q)}.$$
It follows that $\phi(z)=\la\cdot\kap_{-1}(-z,q)$, i.e.
$$c_{12}(z)c_{21}(z)=\frac{c}{2}+\la\cdot\kap_{-1}(-z,q).$$
Substituting $z=q^{-1/2}$ (zero of $c_{21}$) we obtain
$$\frac{c}{2}=-\la\cdot\kap_{-1}(-q^{-1/2},q)=\la\cdot
\kap(-1,-q^{1/2},q)=\la\cdot\frac{\th(1,q)\th(-1,q)}{2}$$
(in the last equality we used (\ref{sp3})).
Similarly, we find
$$c_{11}(z)c_{21}(qz)=\frac{c}{2}-\la\cdot\kap_{-1}(-z,q).$$
So finally the matrix $C$ is equal to
$$\left(\matrix \la\cdot\frac{\th(1,q)\th(-1,q)/2-
\kap_{-1}(-z,q)}{\th(-q^{1/2}z,q^2)} &
\la\cdot\frac{\th(1,q)\th(-1,q)/2+
\kap_{-1}(-z,q)}{\th(-q^{-1/2}z,q^2)}\\ \th(-q^{-1/2}z,q^2)
& -\th(-q^{1/2}z,q^2) \endmatrix \right).$$

Note that since $\det C$ is a non-zero constant as a byproduct we can
find explicitly holomorphic functions $\phi_1$, $\phi_2$ on $\C^*$
with the property
$$\phi_1(z)\th(z,q^2)-\phi_2(z)\th(qz,q^2)=1.$$

\section{Modular property}\label{modular}
 
Throughout this section we use additive notation.
 
\subsection{} Let us recall the  modular interpretation  of
the functional equation for theta-function. For this we have
to consider the stack $\A_1^+$ of elliptic curves (over $\C$) equipped with
a non-trivial point of order 2. It is well-known that
$\A_1^+$ is the quotient of the upper-half plane $\H$ by the
action of the subgroup $\Ga_{1,2}\sub\SL_2(\Z)$ consisting of
matrices
$\ga=\left(\matrix a & b \\ c & d \endmatrix\right)$
with $\det\ga=1$ and $ac\equiv  bd\equiv  0\mod(2)$. The universal elliptic
curve $\EE\ra\A_1^+$ is the quotient of $\C\times\H$ by the
natural  action  of  the  semi-direct product of $\Z^2$ with
$\Ga_{1,2}$, where the action of $\Z^2$ is given by
$$(x,\tau)\mapsto (x+m+n\tau,\tau)$$
while the above matrix in $\Ga_{1,2}$ acts by
$$\ga:(x,\tau)\mapsto (\frac{x}{c\tau+d},
\ga(\tau)=\frac{a\tau+b}{c\tau+d}).$$
Let  $\pi:\C\times\H\ra\EE$ be the natural projection. Below
we will deal with holomorphic vector bundles $V$ on $\EE$ (of ranks 1 and 2)
such  that $\pi^*V$ is trivial. Fixing such a trivialization
gives a 1-cocycle of $G=\Z^2\rtimes\Ga_{1,2}$ with values in the group
of  holomorphic  functions  $\C\times\H\ra\GL_r(\C)$,  where
the action of $G$ on $\C\times\H$ is as above, $r=\rk V$.
Conversely,  given  such  a  1-cocycle $c(g)(x,\tau)$, where
$g\in G$ we can construct a vector  bundle  on  $\EE$  as the
quotient of $\C\times\H\times\C^r$ by the action of $G$ which
sends $((x,\tau),v)$ to $(g(x,\tau),c(g)(x,\tau)v)$.
Thus, the set of isomorphism classes of  holomorphic
bundles on $\EE$ with trivial pull-back to $\C\times\H$ can be
identified with $H^1(G,\GL_r(\O(\C\times\H)))$. Since $G$ is
a semi-direct product of $\Z^2$ and $\Ga_{1,2}$,
any 1-cocycle of $G$ is uniquely determined by its restrictions to
$\Z^2$ and $\Ga_{1,2}$ (these restricted cocycles should be
compatible via the action of $\Ga_{1,2}$ on $\Z^2$).
For  example,  $\th(x,\tau)$ is a global section of the line
bundle $\LL$
on $\EE$ given by a cocycle whose restriction to $\Z^2$ is
$$(m,n)\mapsto \exp(-\pi i n^2\tau-2\pi i n x),$$
and the restriction to $\Ga_{1,2}$ is
\begin{equation}\label{coc-th}
\ga=\left(\matrix a & b \\ c & d \endmatrix\right)\mapsto
\zeta(\ga)\cdot (c\tau+d)^{1/2}\exp(\pi i\frac{cx^2}{c\tau+d})
\end{equation}
where $\zeta(\ga)$ is the 8-th root of unity appearing in the
functional equation for $\th$ (the latter equation expresses
the fact that $\th$ is a section of $\LL$).
Note that the cocycle (\ref{coc-th}) factors into a product
of two cocycles of $\Ga_{1,2}$, the first is with values in
$\O^*(\H)\sub\O^*(\C\times\H)$:
$$\ga\mapsto c_{\th}(\ga)=\zeta(\ga)\cdot(c\tau+d)^{1/2},$$
and the second is given by the exponential factor.
Actually, $c_{\th}$ is a coboundary since $\th(0,\tau)$ is
an invertible function on $\H$ and
$$c_{\th}(\ga)=\frac{\th(0,\ga(\tau))}{\th(0,\tau)}.$$
We  can't  factor  $c_{\th}$  further  due  to
indeterminacy in a choice of a square root of $(c\tau+d)$.
However, $c_{\th}^2$ does factor as the product of a
character $\ga\mapsto \zeta(\ga)^2$ of $\Ga_{1,2}$ with values
in roots of unity of order $4$, and a cocycle
$\ga\mapsto (c\tau+d)$ which corresponds to the line bundle
$\om^{-1}$ on $\A_1^+$ (where $\om$ is the restriction of the
relative canonical bundle $\om_{\EE/\A_1^+}$  to  the  zero
section). In particular, since $c_{\th}^2$ is a coboundary
we obtain that $\om$ is isomorphic to  the  line  bundle  on
$\A_1^+$ associated
with the character $\zeta(\ga)^2$. The explicit formula for
$\zeta(\ga)$ (see e.g. \cite{M}) implies that
$$\zeta(\ga)^2=\cases(-1)^{(d-1)/2}, & \text{$d$ odd},\\
\exp(-\frac{\pi i c}{2}), & \text{$c$ odd}\endcases
$$
 
\subsection{} Now we are going to consider a rank-2 bundle on
$\EE$ which is an extension of $\LL$ by certain line bundle whose
restriction   to  every elliptic  curve  is  of
2-torsion. More presicely, let $L$ be a restriction of $\LL$
to a particular elliptic curve $E$.
Then we have $L\simeq\O(\eta)$  where  $\eta\in  E$  is  the
point   of  order  2  corresponding  to  $(\tau+1)/2$  (this
isomorphism is induced by the theta-function). Now
let $P_{\eta}=\O(\eta-e)\otimes\om^{-1}_E$ be the 2-torsion
line bundle corresponding to  $\eta$,  trivialized  at  zero
$e\in E$. Then we have
$$H^1(L^{-1}\otimes P_{\eta})\simeq
H^1(\O(-e)\otimes\om^{-1}_E)\simeq H^0(\om^2_E)^*.$$
It follows that there is a canonical non-splitting extension
$$0\ra P_{\eta}\otimes\om^2_E\ra V\ra L\ra 0.$$
In other words, there should exist a universal extension on $\EE$
$$0\ra \PP_{\eta}\otimes\om^2\ra\VV\ra\LL\ra 0$$
where $\PP_{\eta}$ is defined in the same way as $P_{\eta}$
with $\eta$ being the divisor $x=(\tau+1)/2\mod(\Z+\Z\tau)$,
$e$ being the zero section, and $\om_{E}$ replaced by $\om$.
The pull-back of $\PP_{\eta}$ to $\C\times\H$ is trivial and
the corresponding 1-cocycle of $G$ is given by
$$(m,n)\mapsto \exp(\pi i n(\tau+1)),$$
$$\left(\matrix a & b \\ c & d \endmatrix\right)\mapsto
\exp(\pi i(\frac{1}{c\tau+d}-1)x).$$
Below we construct a bundle $\VV$ on $\EE$ fitting into above
exact sequence together
with a trivialization of $\pi^*\VV$ such that the pair
$(\kap,\th)$ defines a global section of $\VV$.
Note  that  one  has a similar rank-2 bundle on the stack of
elliptic curves (without a choice of a point of order 2),
which is an extension of $\O(e)$ by $\om$. However, we will
stick to the case of $\A_1^+$ in order to stress the
analogy with the usual functional equation for $\th$
(which takes place on this stack).
 
\subsection{}
Recall that for an elliptic curve $E=E^{\tau}=\C/(\Z+\Z\tau)$
we have defined in section \ref{rank2} the rank-2 bundle
on $E^{\tau}$
$$V^{\tau}=F_{\exp(\pi i(\tau+1))}=V_2
\left(  \matrix \exp(\pi i(\tau+1)) & 1 \\ 0 &
\exp(-\pi i \tau-2\pi i x) \endmatrix\right)$$
which is an extension of $L=L^{\tau}$ by
$P_{\eta}=P_{\eta}^{\tau}$. Now
for an element
$$\ga=\left(\matrix a & b \\ c & d \endmatrix\right)
\in\Ga_{1,2}$$
we have the corresponding isomorphism
$$f_{\ga}:E^{\tau}\ra E^{\ga(\tau)}:x\mapsto \frac{x}
{c\tau+d}$$
induced by the action of $\Ga_{1,2}$ on $\C\times\H$.
The 1-cocycle of $G$ corresponding to  $\LL$ (resp.
$\PP_{\eta}$) evaluated at $\ga$ induces an isomorphism
$f_{\ga}^*L^{\ga(\tau)}\simeq L^{\tau}$ (resp.
$f_{\ga}^*P_{\eta}^{\ga(\tau)}\simeq P_{\eta}^{\tau}$).
Using these isomorphisms we consider
the pull-back by $f_{\ga}$ of the extension
$$0\ra P_{\eta}^{\ga(\tau)}\ra V^{\ga(\tau)}\ra L^{\ga(\tau)}
\ra 0$$
as an extension of $L^{\tau}$ by $P_{\eta}^{\tau}$. The class of the
latter extension differs from the class of the extension given
by $V^{\tau}$ by a non-zero constant. Hence, there
exists a unique isomorphism
$$V^{\tau}\wt{\ra} f_{\ga}^*V^{\ga(\tau)}$$
which   induces   the  identity  map  on  $L^{\tau}$  and  a
a  non-zero   constant   multiple   of   the   identity   on
$P_{\eta}^{\tau}$.
 
Equivalently, one can say that there exists a unique 1-cocycle
of $G$ with values in $\GL_2(\O(\C\times\H))$ such that its restriction to
$\Z^2$ is
\begin{equation}
(m,n)\mapsto
\left(\matrix \exp(\pi i n(\tau+1)) & 1 \\ 0 &
\exp(-\pi i n^2\tau-2\pi i n x) \endmatrix\right)
\end{equation}
and its restriction to $\Ga_{1,2}$ has form
\begin{equation}
\ga=\left(\matrix a & b \\ c & d \endmatrix\right)
\mapsto
\left(\matrix k_{\ga}(\tau)\cdot\exp(\pi i(\frac{1}{c\tau+d}-1)x)
& \phi_{\ga}(x,\tau) \\ 0 &
c_{\th}(\ga)\cdot\exp(\pi i \frac{cx^2}{c\tau+d}) \endmatrix\right)
\end{equation}
for  some holomorphic function $\phi_{\ga}(x,\tau)$
and some invertible holomorphic function $k_{\ga}(\tau)$.
Moreover, it follows from definition of $\kap$ that
$$\left[\matrix \kap((\tau+1)/2,x,\tau) \\ \th(x,\tau)
\endmatrix\right]$$
is a section of the bundle $\VV$ defined by the above cocycle.
This allows us to find functions $k_{\ga}$ and $\phi_{\ga}$
explicitly. Namely, for any $\ga\in\Ga_{1,2}$ we should have
\begin{equation}\label{div}
\kap((\ga(\tau)+1)/2,\frac{x}{c\tau+d},\ga(\tau))=
k_{\ga}(\tau)\cdot\exp(\pi i(\frac{1}{c\tau+d}-1)x)\cdot
\kap((\tau+1)/2,x,\tau)+\phi_{\ga}(x,\tau)\th(x,\tau).
\end{equation}
To find $k_{\ga}$ let us substitute $x=(\tau+1)/2$ in this equation.
Note that the equation (\ref{def}) together with the vanishing
$\th((\tau+1)/2,\tau)=0$ implies that
$$\kap((\tau+1)/2,x+m+n\tau,\tau)=\exp(\pi i n(\tau+1))
\kap((\tau+1)/2,x,\tau).$$
Using this together with the formulas
$$\frac{1}{c\tau+d}=a-c\ga(\tau),$$
$$\frac{\tau}{c\tau+d}=-b+d\ga(\tau)$$
we get from (\ref{div})
\begin{align*}
&\exp(\pi i\frac{d-c-1}{2}(\ga(\tau)+1))
\kap((\ga(\tau)+1)/2,(\ga(\tau)+1)/2,\ga(\tau))=\\
&k_{\ga}(\tau)\cdot\exp(\frac{\pi i}{2}(a-b-1+
(d-c)\ga(\tau)-\tau))\cdot
\kap((\tau+1)/2,(\tau+1)/2,\tau),
\end{align*}
i.e.
$$k_{\ga}=\exp(\frac{\pi i}{2}(d-a+b-c-\ga(\tau)+\tau))\cdot
\frac{\kap((\ga(\tau)+1)/2,(\ga(\tau)+1)/2,\ga(\tau))}
{\kap((\tau+1)/2,(\tau+1)/2,\tau)}.$$
Now using the formula (\ref{sp2}) for
$\kap((\tau+1)/2,(\tau+1)/2,\tau)$ and the
functional equation for the theta-function we obtain
$$k_{\ga}=\exp(\frac{\pi i}{2}(d-a+b-c)+
\frac{3\pi i}{4}(\tau-\ga(\tau)))\cdot
\zeta(\ga)^2\cdot\chi(\ga)\cdot (c\tau+d)$$
where $\chi$ is the character of $\Ga_{1,2}$ with
values in 4th roots of unity defined as follows:
$$\chi(\ga)=\cases (-1)^{a/2}\exp(\frac{\pi i}{4}
(ab+cd)), & \text {$a$ even},\\
(-1)^{c/2}\exp(\frac{\pi i}{4}(ab+cd)), & \text{$c$ even}
\endcases$$
It is easy to see that $(-1)^{(d-a+b-c)/2}=\chi(\ga)^2\zeta(\ga)^4$,
so we can write finally
\begin{equation}\label{kga}
k_{\ga}=\exp(\frac{3\pi i}{4}(\tau-\ga(\tau)))\cdot
\zeta(\ga)^{-2}\cdot\chi^{-1}(\ga)\cdot (c\tau+d)
\end{equation}
Now the equation (\ref{div}) gives the explicit
formula for $\phi_{\ga}$.
 
Summarizing we can say
that the existence of the rank-2 bundle $\VV$ on the
universal elliptic curve $\EE$ is equivalent to the following
divisibility property of $\kap$. Let us denote
$$\kap_0(x,\tau)=\exp(\frac{3\pi i\tau}{4})
\kap((\tau+1)/2,x,\tau).$$
Then for every $\ga\in\Ga_{1,2}$ the difference
$$\kap_0(\frac{x}{c\tau+d},\ga(\tau))-
\zeta(\ga)^{-2}\cdot\chi^{-1}(\ga)\cdot (c\tau+d)
\cdot\exp(\pi i(\frac{1}{c\tau+d}-1)x)\cdot
\kap_0(x,\tau)$$
is divisible by $\th(x,\tau)$.

\end{document}